\documentclass[12pt]{amsart}

\usepackage{amssymb}
\usepackage{amscd}
\swapnumbers
\def\frak{\mathfrak}
\def\Bbb{\mathbb}
\def\Cal{\mathcal}

\newtheorem*{prop*}{Proposition \thesubsection}

\newtheorem*{thm*}{Theorem \thesubsection}

\newtheorem*{lem*}{Lemma}

\newtheorem*{kor*}{Corollary \thesubsection}

\newcommand{\gr}{\operatorname{gr}}

\newcommand{\Ad}{\operatorname{Ad}}

\newcommand{\im}{\operatorname{im}}

\newcommand{\Aut}{\operatorname{Aut}}

\newcommand{\fg}{{\frak g}}
\newcommand{\x}{\times}
\renewcommand{\o}{\circ}

\let\ccdot\cdot
\def\cdot{\hbox to 2.5pt{\hss$\ccdot$\hss}}

\newcommand{\g}{{\frak g}}

\newcommand{\al}{\alpha}
\newcommand{\be}{\beta}

\newcommand{\ka}{\kappa}

\newcommand{\om}{\omega}
\renewcommand{\phi}{\varphi}
\newcommand{\ph}{\varphi}
\newcommand{\ps}{\psi}

\newcommand{\Ga}{\Gamma}
\newcommand{\La}{\Lambda}
\newcommand{\Om}{\Omega}
\newcommand{\Ph}{\Phi}

\def\sideremark#1{\ifvmode\leavevmode\fi\vadjust{
\vbox to0pt{\hbox to 0pt{\hskip\hsize\hskip1em
\vbox{\hsize3cm\tiny\raggedright\pretolerance10000
\noindent #1\hfill}\hss}\vbox to8pt{\vfil}\vss}}}

\begin{document}
\title[Automorphism group of distributions]{On automorphism groups of
  some types of generic distributions} 

\author{Andreas \v Cap and Katharina Neusser}

\thanks{First author supported by project P19500-N13 of the
    ``Fonds zur F\"orderung der wissenschaftlichen Forschung'' (FWF),
    second author supported by Initiativkolleg IK-1008 of the
    University of Vienna.}

\address{Fakult\"at f\"ur Mathematik, Universit\"at Wien, Nordbergstr\ss e 15, A--1090 Wien, Austria}

\email{Andreas.Cap@esi.ac.at, Katharina.Neusser@univie.ac.at}



\begin{abstract}
  To certain types of generic distributions (subbundles in a tangent
  bundle) one can associate canonical Cartan connections. Many of
  these constructions fall into the class of parabolic geometries. The
  aim of this article is to show how strong restrictions on the
  possibles sizes of automorphism groups of such distributions can be
  deduced from the existence of canonical Cartan connections.
  This needs no information on how the Cartan connections
  are actually constructed and only very basic information on their
  properties. In particular, we discuss the examples of generic
  distributions of rank two in dimension five, rank three in dimension
  six, and rank four in dimension seven. 
\end{abstract}

\maketitle

\section{Introduction}\label{1}

This article deals with geometric questions on subbundles in the
tangent bundles of smooth manifolds. While such structures always have
been of interest in control theory, their importance in various parts
of differential geometry and geometric analysis has increased a lot
during the last years. This refers to, for example, sub--Riemannian
structures, questions related to Carnot groups and
Carnot--Caratheodory manifolds, as well as analytical properties of
differential operators obtained as sums of squares of sections of such
subbundles.

Since integrable subbundles can be removed by passing to leaves of the
corresponding foliation, one usually restricts the attention to
\textit{bracket generating} distributions. This condition means that
sections of the distribution together with their iterated Lie brackets
span the full tangent bundle.

Even under this assumption, different types of distributions can have
entirely different behavior. For example, consider automorphisms of a
distribution, i.e.~diffeomorphisms of the manifold whose derivatives
in all points respect the distribution. On the one hand, there are
examples like contact distributions, which admit a local normal form
and always have infinite dimensional families of automorphisms. The
first example of the other possible behavior was found independently
by F.~Engel and E.~Cartan in their work on exceptional Lie algebras of
type $G_2$. It was studied in detail in Cartan's famous ``five
variables paper'' \cite{Cartan:five}. In this article, he studied
distributions of rank two and three on five dimensional manifolds,
which, in addition to being bracket generating, satisfy a genericity
condition. He associated to such distributions a canonical
\textit{Cartan connection} on a certain principal bundle.  This
immediately implies that such distributions have local invariants
(similar to the curvature of a Riemannian metric) and hence cannot
admit simple local normal forms. Further, it implies that the
automorphisms of such a distribution form a finite dimensional Lie
group and each automorphism is determined by some finite jet (actually
the two--jet) in one point.

Cartan's result has been (much later) extended to various other
generic types of distributions. Many of these examples fall into the
class of so called \textit{parabolic geometries}, since the homogeneous
model of the geometry is the quotient of a semisimple Lie group by a
parabolic subgroup. Motivated by the examples of conformal structures
and CR structures, these geometries have been intensively studied
during the last years, and many striking results have been achieved. 

The results on existence of canonical Cartan connections are usually
difficult and Cartan connections themselves are often considered as
being hard to use. In this article we want to show that Cartan
connections lead to interesting results on distributions in a rather
simple way. For these applications, no knowledge about the actual
construction of the Cartan connections but only some basic information
about their properties is needed. We show that simple algebraic
computations (mostly linear algebra) can be used to obtain surprising
restrictions on the possible dimensions of the automorphism groups of
certain types of distributions.

The basic ideas we use certainly go back to Cartan, the more specific
version for parabolic subalgebras has been implicitly used in
\cite{Yamaguchi:CR} in the study of automorphism groups of CR
structures. They have been explicitly formulated in \cite{Srni04} in
the context of general parabolic geometries. The algebraic
computations needed to apply these ideas in the cases of generic
distributions discussed here as well as some of the realizations of
automorphism groups are part of the second author's diploma thesis,
see \cite{Neusser}.

\section{Distributions which are equivalent to parabolic
  geometries}\label{2} In this section, we give a brief description of
the relation between certain types of distributions and parabolic
geometries, i.e.~Cartan connections with homogeneous model the
quotient of a semisimple Lie group by a parabolic subgroup. 

\subsection{Parabolic subalgebras}\label{2.1}
These are a special type of subalgebras in semisimple Lie algebras,
which can be defined in several equivalent ways. In terms of structure
theory, one best defines a parabolic subalgebra in a complex
semisimple Lie algebra as one which contains a maximal solvable
subalgebra (which usually is called a \textit{Borel subalgebra}). Then
one defines parabolic subalgebras in real semisimple Lie algebras via
complexification. Equivalently, one may define a subalgebra in an
arbitrary semisimple Lie algebra to be parabolic if its nilradical
coincides with its annihilator under the Killing form, see \cite{CDS}. 

For our purposes, the most useful definition is the one in terms of
$|k|$--gradings, which also handles the real and complex cases
simultaneously. 

\subsection*{Definition}
Let $\fg$ be a real or complex semisimple Lie algebra, and let $k$ be
a positive integer. 

\noindent
(1) A $|k|$--grading on $\fg$ is a vector space decomposition
$\fg=\fg_{-k}\oplus\dots\oplus\fg_k$ such that we have
$[\fg_i,\fg_j]\subset\fg_{i+j}$ for all $i$ and $j$, with the
convention that $\fg_\ell=\{0\}$ for $|\ell|>k$, and such that the
subalgebra $\fg_-:=\fg_{-k}\oplus\dots\oplus\fg_{-1}$ is generated by
$\fg_{-1}$. 

\noindent
(2) Given a $|k|$--grading as in (1), we put
    $\fg^i:=\fg_i\oplus\dots\oplus\fg_k$ for all $i$ as well
    as $\frak p:=\fg^0$ and $\frak p_+:=\fg^1$. 

\noindent
(3) A Lie subalgebra $\frak p$ of $\fg$ is called \textit{parabolic}
if it can be realized as $\fg^0$ for some $|k|$--grading of $\fg$.

\medskip

The subspaces $\fg^i$ from (2) define a decreasing filtration of
$\fg$, which makes $\fg$ into a \textit{filtered Lie algebra},
i.e.~$[\fg^i,\fg^j]\subset\fg^{i+j}$ for all $i,j$. In particular,
this implies that $\frak p=\fg^0$ is a Lie subalgebra of $\fg$, and
that $\frak p_+=\fg^1$ is an ideal in $\frak p$, which is nilpotent
since the grading has finite length. Likewise, $\fg_-\subset\fg$ is a
nilpotent Lie subalgebra. It turns out, see \cite{Yamaguchi}, that the
Lie algebras $\fg_-$ and $\frak p_+$ are always isomorphic. From the
filtration property it also follows that each of the filtration
components $\fg^i$ is $\frak p$--invariant.

By the grading property, $\fg_0\subset\frak p\subset\fg$ is a Lie
subalgebra, and each of the grading components $\fg_i$ is
$\fg_0$--invariant. It turns out that the Lie algebra $\fg_0$ is
always \textit{reductive}, so it is the direct sum of a semisimple Lie
algebra and a center. In particular, the representation theory of
$\fg_0$ is significantly easier than the one of the parabolic $\frak
p$, which makes $\fg_0$ a valuable technical tool in the
theory. Moreover, it turns out that completely reducible
representations of $\frak p$ always come from representation of
$\fg_0$ via the quotient map $\frak p\to \frak p/\frak
p_+\cong\fg_0$. 

A crucial feature of the theory is the computability of the Lie
algebra cohomology groups $H^*(\fg_-,\fg)$ via Kostant's version of
the Bott--Borel--Weil theorem, see \cite{Kostant}. The standard
complex for computing these cohomologies consists of spaces of
multilinear alternating maps from $\fg_-$ to $\fg$. In particular,
$\fg_0$ acts naturally on these spaces and it is easy to see that the
differentials in the standard complex are $\fg_0$--equivariant. Hence
the cohomologies naturally are representations of $\fg_0$, and
Kostant's results describes them as such representations. The
computation of the cohomologies is completely algorithmic and it is
also implemented by J.~\v Silhan as an extension to the Lie software
system, see \cite{Silhan}, so the computation of the cohomologies can
be left to a computer. 

\subsection{The symbol algebra of a distribution}\label{2.2}
Let $M$ be a smooth manifold and let $H\subset TM$ be a distribution,
i.e.~a smooth subbundle. We assume that $H$ is bracket generating,
i.e.~that sections of $H$ together with their iterated Lie brackets
span the whole tangent bundle. We will also write $T^{-1}M$ for
$H$. Next we require that sections of $H$ together with Lie brackets
of two such sections span a smooth subbundle $T^{-2}M\subset TM$,
which by construction contains $T^{-1}M$. Inductively, we require that
we get a filtration $TM=T^{-k}M\supset T^{-k+1}M\supset\dots\supset
T^{-1}M$ of the tangent bundle by smooth subbundles, such that for
each $i<0$ sections of $T^iM$ together with Lie brackets of one
section of $T^iM$ and one section of $T^{-1}M$ span $T^{i-1}M$. The
sequence of the ranks of the subbundles $T^iM$ is usually called the
\textit{small growth vector} of the distribution $H$. 

Having extended the distribution $H$ to the filtration $\{T^iM\}$, one
can next encode the non--integrability properties of $H$. Namely, for
each $i=-k,\dots,-1$, one defines $\gr_i(TM):=T^iM/T^{i+1}M$, and then
$\gr(TM)=\oplus\gr_i(TM)$ is the associated graded to the filtered
vector bundle $TM$. One immediately verifies, that for sections
$\xi\in\Ga(T^iM)$ and $\eta\in\Ga(T^jM)$, the Lie bracket $[\xi,\eta]$
is a section of $T^{i+j}M$, where $T^\ell M=TM$ for $\ell\leq
-k$. Hence the Lie bracket of vector fields induces a bilinear bundle
map $\gr_i(TM)\x\gr_j(TM)\to\gr_{i+j}(TM)$, which for each $x\in M$
makes $\gr(T_xM)$ into a nilpotent graded Lie algebra. This is called
the \textit{symbol algebra} of the distribution $H$ at $x$. 

Suppose that $f:M\to M$ is a diffeomorphism such that $Tf(H)\subset
H$. Then by construction $Tf$ preserves each of the subbundles $T^iM$
and hence is compatible with the filtration of $TM$. Hence for each
$x\in M$, the tangent map $T_xf$ induces a linear isomorphism
$\gr(T_xM)\to\gr(T_{f(x)}M)$. Compatibility of $Tf$ with the Lie
bracket of vector fields immediately implies that this map actually is
an isomorphism of the symbol algebras at $x$ and $f(x)$. Hence the
symbol algebra is a fundamental invariant of the distribution.

In general, the isomorphism type of the symbol algebra may change from
point to point, but we will be only interested in distributions for
which $\gr(TM)$ is locally trivial as a bundle of nilpotent graded Lie
algebras. If $\frak n=\frak n_{-k}\oplus\dots\oplus\frak n_{-1}$ is
the modelling nilpotent graded Lie algebra, then $\gr(TM)$ has an
obvious natural frame bundle with structure group $\Aut_{gr}(\frak
n)$, the group of automorphisms of the graded Lie algebra $\frak n$,
compare with \cite{Morimoto:filtered}.

\subsection{From parabolics to distributions}\label{2.3}
Consider a $|k|$--graded Lie algebra $\fg=\oplus_{i=-k}^k\fg_i$ as in
\ref{2.1} and the corresponding filtration $\{\fg^i\}$. Let $G$ be a
Lie group with Lie algebra $\frak g$. Then one shows, see
\cite{Yamaguchi} and \cite{C-S}, that 
\begin{gather*}
  P:=\{g\in G:\Ad(g)(\fg^i)\subset\fg^i\quad\forall i\}\\
  G_0:=\{g\in G:\Ad(g)(\fg_i)\subset\fg_i\quad\forall i\}
\end{gather*}
are closed subgroups of $G$ with Lie algebras $\frak p=\fg^0$ and
$\fg_0$, respectively. In particular, restricting the action of $G_0$
to the nilpotent subalgebra $\fg_-$ we obtain an action by Lie algebra
automorphisms, i.e.~a homomorphism $G_0\to\Aut_{gr}(\fg_-)$. 

The homogeneous spaces of the form $G/P$ are the so--called
\textit{generalized flag manifolds}, which are of central interest in
representation theory. It turns out that they are always compact.

\begin{prop*}
  The generalized flag manifold $G/P$ carries a natural distribution
  $H\subset T(G/P)$ of rank $\dim(\fg_{-1})$, whose bundle of symbol
  algebras is locally trivial with modelling algebra $\fg_-$. The
  natural left action of $G$ on $G/P$ preserves this distribution. 
\end{prop*}
\begin{proof}
  The tangent bundle of $G/P$ can be identified with the associated
  bundle $G\x_P(\fg/\frak p)$. Now the $P$--invariant filtration
  $\{\fg^i\}$ of $\fg$ induces a $P$--invariant filtration $\fg/\frak
  p=\fg^{-k}/\frak p\supset\dots\supset\fg^{-1}/\frak p$. For each
  $i=-k,\dots,-1$, we get a smooth subbundle $G\x_P(\fg^i/\frak
  p)=:T^i(G/P)\subset T(G/P)$,~i.e. a filtration of the tangent bundle
  of $G/P$. Explicitly, denoting by $p:G\to G/P$ the natural
  projection, the subbundle $T^i(G/P)$ is given by
$$
T_{gP}^i(G/P)=T_gp(\{L_X(g):X\in\fg^i\}),
$$
where $L_X$ denotes the left invariant vector field generated by
$X$. This immediately shows that the left actions of elements of $G$
preserve the filtration $\{T^i(G/P)\}$ of the tangent bundle. 

By construction, the component $\gr_i(T(G/P))$ of the associated
graded is the bundle induced by the representation $(\fg^i/\frak
p)/(\fg^{i+1}/\frak p)\cong \fg^i/\fg^{i+1}$. Consider sections
$\xi\in T^i(G/P)$ and $\eta\in T^j(G/P)$ with $i<j$. By construction,
there are local lifts $\tilde\xi,\tilde\eta\in\frak X(G)$ which can be
written in the form $\tilde\xi=\sum_a\phi_aL_{X_a}$ and
$\tilde\eta=\sum_b\ps_bL_{Y_b}$ for smooth functions $\phi_a$ and
$\ps_b$ and elements $X_a\in\fg^i$ and $Y_b\in\fg^j$. This shows that
$[\tilde\xi,\tilde\eta]$ can be written as the sum of 
$$
\textstyle\sum_{a,b}\ph_a\ps_b[L_{X_a},L_{Y_b}]=
\sum_{a,b}\ph_a\ps_bL_{[X_a,Y_b]} 
$$
and a linear combination of left invariant vector fields with
generators in $\fg^i$. Hence we conclude that
$[\xi,\eta]\in\Ga(T^{i+j}(G/P))$. Since $\fg_-$ is generated by
$\fg_{-1}$ we conclude that the distribution $H:=T^{-1}(G/P)$ is
bracket generating and that $\{T^i(G/P)\}$ is the associated
filtration as described in \ref{2.2}. Finally, it also shows that
under the natural identification $\fg^i/\fg^{i+1}\cong\fg_i$, the
symbol algebra of $H$ in each point is isomorphic to $\fg_-$.
\end{proof}

\subsection{Canonical Cartan connections}\label{2.4}
Since the algebras $\fg$ are always semisimple, there is a natural
choice of a Lie group $G$ with Lie algebra $\fg$, namely the
automorphism group $\Aut(\fg)$. Then the subgroups $G_0\subset
P\subset G$ are the groups $\Aut_{gr}(\fg)\subset\Aut_f(\fg)$ of
automorphisms preserving the grading respectively the filtration of
$\fg$. Now we have to assume an additional (cohomological) condition
on the grading of $\fg$. The first Lie algebra cohomology group
$H^1(\fg_-,\fg)$ consists of equivalence classes of linear maps
$\fg_-\to\fg$, and there is an obvious notion of homogeneity for such
maps. The Lie algebra differentials are compatible with homogeneity,
so the cohomology group naturally splits according to homogeneity.

Assuming that $H^1(\fg_-,\fg)$ is concentrated in negative
homogeneities, it turns out (see \cite{katja}) that the homomorphism
$G_0\to\Aut_{gr}(\fg_-)$ is actually an isomorphism, and then the
general prolongation procedures of \cite{Tanaka,Morimoto,C-S} imply 
\begin{thm*}
  Let $M$ be a smooth manifold such that $\dim(M)=\dim(G/P)$ endowed
  with a bracket generating distribution $H\subset TM$ of rank
  $\dim(\fg_{-1})$, whose bundle of symbol algebras is locally trivial
  with modelling Lie algebra $\fg_-$. Then the natural frame bundle
  for $M$ with structure group $\Aut_{gr}(\fg_-)\cong G_0$ can be
  canonically extended to a principal $P$--bundle $\Cal G\to M$, which
  can be endowed with a regular normal Cartan connection
  $\om\in\Om^1(\Cal G,\fg)$.
  
  The pair $(\Cal G,\om)$ is uniquely determined up to isomorphism,
  and the construction actually establishes an equivalence of
  categories between manifolds endowed with appropriate distributions
  and regular normal Cartan geometries.
\end{thm*}

Let us explain this a bit. First of all, it is easy to see (compare
with \cite{Yamaguchi,C-S}) that $G_0$ can be naturally viewed as a
quotient of $P$. Indeed, the exponential map defines a diffeomorphism
from $\fg^1$ onto a closed subgroup $P_+\subset P$ such that
$P/P_+\cong G_0$. The statement that $\Cal G$ extends the natural
frame bundle then simply means that the quotient $\Cal G/P_+$ is (via
forms induced by the Cartan connection) isomorphic to this frame
bundle.

A Cartan connection on $\Cal G$ by definition is a one--form
$\om\in\Om^1(\Cal G,\fg)$ such that $\om(u):T_u\Cal G\to\fg$ is a
linear isomorphism for each $u\in\Cal G$. Further, $\om$ has to be
equivariant for the principal right action of $P$,
i.e.~$(r^g)^*\om=\Ad(g^{-1})\o\om$ for all $g\in P$, and it has to
reproduce the generators of fundamental vector fields. The pair $(\Cal
G,\om)$ is then referred to as a \textit{Cartan geometry} of type
$(G,P)$.  Cartan geometries can be viewed as ``curved analogs'' of the
homogeneous space $G/P$, which determines a Cartan geometry via the
natural principal bundle $G\to G/P$ and the left Maurer--Cartan form
on $G$. In this context, $G/P$ is referred to as the
\textit{homogeneous model} of geometries of type $(G,P)$.

The amount to which a general Cartan geometry differs from the
homogeneous model is measured by its \textit{curvature}. This is the
two--form $K\in\Om^2(\Cal G,\fg)$ defined by
$$
K(\xi,\eta)=d\om(\xi,\eta)+[\om(\xi),\om(\eta)]. 
$$
From the defining properties of a Cartan connection it immediately
follows that $K$ is horizontal and $P$--equivariant. In particular, its
value on $\xi$ and $\eta$ depends only on the projections of the
tangent vectors to $M$. The uniqueness of the Cartan connection in the
theorem is ensured by the conditions of regularity and normality on the
curvature. Regularity means that if the image of $\xi$ and $\eta$ in
$TM$ lie in the subbundles $T^iM$ and $T^jM$, respectively, then
$K(\xi,\eta)\in \fg^{i+j+1}$. 

The condition of normality is crucial for the uniqueness question, and
finding appropriate normalization conditions often is a very difficult
step in the construction of canonical Cartan connections. For the
purposes of this article, we do not need any details on the form of
this condition. The only information we need (and also this is only
needed to deal with the non--flat case) is that $K$ determines a
natural quantity called the \textit{harmonic curvature}, which is a
section $\ka_h$ of the associated bundle $\Cal G\x_P(H^2(\fg_-,\fg))$.
Vanishing of the harmonic curvature is equivalent to vanishing of $K$
and to local isomorphism of the geometry with $G/P$, see \cite{C-S}. 

The theorem on existence and uniqueness of Cartan connections has
important consequences for the homogeneous model.
\begin{kor*}
  Suppose that $H^1(\fg_-,\fg)$ is concentrated in negative
  homogeneities. Then for any connected open subset $U\subset G/P$,
  the automorphism group of the distribution
  $T^{-1}U:=T^{-1}(G/P)|_U\subset TU$ is the subgroup of $G$ consisting
  of all elements whose left action on $G/P$ maps the subset $U$ to
  itself. In particular, the automorphism group of $T^{-1}(G/P)$
  itself is the group $G=\Aut(\fg)$.
\end{kor*}
\begin{proof}
  Of course, any element of $G$ whose left action preserves $U$ gives
  rise to an automorphism. Conversely, let $p:G\to G/P$ be the natural
  projection and let $\om$ be the left Maurer Cartan form on $G$. Then
  $\om$ restricts to a Cartan connection on the principal $P$--bundle
  $p^{-1}(U)\to U$, which is flat by the Maurer--Cartan equation.
  Hence it must be the normal Cartan connection determined by
  $T^{-1}U$. The equivalence of categories stated in the theorem
  implies that any automorphism of $U$ lifts to an automorphism of the
  Cartan geometry, i.e.~to a $P$--equivariant diffeomorphism on
  $p^{-1}(U)$ which is compatible with the Cartan connection. Then 
  \cite[Theorem 5.2]{Sharpe} applied to the given automorphism and the 
  identity shows that on each connected component of $p^{-1}(U)$ the 
  automorphism is given by the left action of some element of $G$. Since $U$ is 
  connected, $P$ acts transitively on the set of connected components of
  $p^{-1}(U)$, so $P$--equivariancy implies that we always get the
  same element $g\in G$.
\end{proof}

\subsection{Automorphism groups}\label{2.5}
For distributions which are equivalent to Cartan geometries, any
automorphism of the distributions canonically lifts to the Cartan
geometry, so the automorphism groups of the distribution can be
identified with the one of the Cartan geometry. Of course, an
automorphism of a Cartan geometry $(\Cal G,\om)$ is a $P$--equivariant
diffeomorphism $\Ph:\Cal G\to\Cal G$ such that $\Ph^*\om=\om$. There
is an obvious infinitesimal analog of this concept. Namely, the flow
of a complete vector field on $\Cal G$ consists of
automorphisms if and only if the field lies in
$$
\frak{inf}(\Cal G,\om):=\{\xi\in\frak X(\Cal G):(r^g)^*\xi=\xi \text{\
  for all\ }g\in P\text{\
  and\ } \Cal L_\xi\om=0\},
$$
where $(r^g)$ is the principal right action of $g$ and $\Cal L$
denotes the Lie derivative. This is called the set of
\textit{infinitesimal automorphisms} of the geometry. Notice that by
definition $\frak{inf}(\Cal G,\om)$ is closed under Lie brackets.

\begin{thm*}
  (1) Let $(p:\Cal G\to M,\om)$ be a Cartan geometry of type $(G,P)$
  with connected base $M$.  Then the automorphism group $\Aut(\Cal
  G,\om)$ can be made into a Lie group of dimension $\leq\dim(G)$,
  whose Lie algebra $\frak{aut}(\Cal G,\om)$ consists of all complete
  vector fields contained in $\frak{inf}(\Cal G,\om)$.

\noindent
(2) For any point $u\in\Cal G$, the map $\xi\mapsto\om(\xi(u))$
    induces an injection $\frak{aut}(\Cal G,\om)\hookrightarrow\frak
    g$. Denoting by $\frak a\subset\frak g$ the image, the Lie bracket
    on $\frak{aut}(\Cal G,\om)$ is mapped to the operation
    $$(X,Y)\mapsto [X,Y]-K(\om^{-1}(X),\om^{-1}(Y))(u)$$
    on $\frak a$. 

\noindent
(3) If the Cartan geometry is regular, then restricting the filtration
$\{\fg^i\}$ of $\fg$ to the subspace $\frak a$ makes $\frak{aut}(\Cal
G,\om)$ into a filtered Lie algebra. The associated graded of this Lie
algebra is isomorphic to a graded Lie subalgebra of $\fg$.
\end{thm*}
\begin{proof}
  We first claim that $\xi\mapsto\omega(\xi(u))$ defines an injection
  $\frak{inf}(\Cal G,\om)\to\fg$, i.e.~that any infinitesimal
  automorphism is determined by its value in one point. Using the
  standard formula for the Lie derivative and inserting the definition
  of the exterior derivative, we get
$$
0=(\Cal L_\xi\om)(\eta)=d\om(\xi,\eta)+\eta\cdot\om(\xi)=
\xi\cdot\om(\eta)-\om([\xi,\eta]).
$$
for $\xi\in\frak{inf}(\Cal G,\om)$ and $\eta\in\frak X(\Cal G)$. In
particular, if $\om(\eta)$ is constant, then injectivity of $\om$
implies $[\xi,\eta]=0$. But this implies that $\xi$ is invariant under
the flows of such fields. Since such fields span each tangent space,
we see that $\xi(u)$ determines $\xi$ locally around $u$. Since $\xi$
is $P$--invariant, this determines the restriction of $\xi$ to
$p^{-1}(U)$ for some open neighborhood $U$ of $p(u)$ in $M$. Now
connectedness of $M$ implies that $\xi$ is globally determined by
$\xi(u)$. In particular, $\frak{inf}(\Cal G,\om)$ is a finite
dimensional Lie subalgebra of $\frak X(\Cal G)$, and (1) follows from
R.~Palais' characterization of Lie transformation groups, see
\cite{Palais}.

\noindent
(2) We already know that $\frak{aut}(\Cal G,\om)$ injects into $\fg$,
    so it remains to verify the formula for the Lie bracket. It is
    well known that the Lie bracket on $\frak{aut}(\Cal G,\om)$ is
    induced by the negative of the Lie bracket of vector fields. Now
    for $\xi,\eta\in\frak{inf}(\Cal G,\om)$, we get 
$$
0=d\om(\xi,\eta)+\eta\cdot\om(\xi)=K(\xi,\eta)-[\om(\xi),\om(\eta)]
+\om([\eta,\xi]), 
$$
where we have used the definition of $K$ to rewrite the first term and
the infinitesimal automorphism equation for $\eta$ to rewrite the
second one. This evidently implies that claimed formula. 

\noindent
(3) The regularity of the geometry implies that for $X\in\fg^i$ and
    $Y\in\fg^j$ we get
    $K(\om^{-1}(X),\om^{-1}(Y))(u)\in\fg^{i+j+1}$. This implies that
    the modified bracket on $\frak a$ still respects the filtration
    and that the $K$--term does not contribute to the bracket on the
    associated graded. Since the filtration of $\fg$ comes from a
    grading, the associated graded of $\fg$ is isomorphic to $\fg$
    itself, and the last statement follows.   
\end{proof}

Part (3) of this theorem is the main input for getting restrictions on
possible automorphism groups. To get additional information for
geometries which are non--flat (i.e.~not locally isomorphic to $G/P$)
we need an additional bit of information. We have briefly discussed in
\ref{2.4} the harmonic curvature $\ka_h$ of $\om$ which is a section
of $\Cal G\x_P H^2(\fg_-,\frak g)$. Hence it corresponds to a
$P$--equivariant function $\Cal G\to H^2(\fg_-,\frak g)$. Now it turns
out that the latter representation is completely reducible, so the
$P$--action factorizes through $G_0$. Naturality of the construction of
$\ka_h$ implies that any infinitesimal automorphism of $(\Cal G,\om)$
has to preserve $\ka_h$, i.e.~annihilate the corresponding
function. Using this, we formulate

\begin{kor*}
  Let $\fg=\oplus_{i=-k}^k\fg_i$ be a $|k|$--graded Lie algebra such
  that $H^1(\fg_-,\fg)$ is concentrated in negative homogeneous
  degrees, and let $G_0\subset P\subset G$ be the corresponding
  groups. Let $M$ be a smooth manifold of dimension $\dim(G/P)$
  endowed with a distribution $H\subset TM$ of rank $\dim(\fg_{-1})$
  such that the bundle of symbol algebras is locally trivial and
  modelled on $\fg_-$. Let $\Aut(M,H)$ be the group of all
  diffeomorphisms of $M$ which preserve the distribution $H$. 

\noindent
(1) $\Aut(M,H)$ is a Lie group of dimension $\leq\dim(G)$. 

\noindent
(2) If $\ell:=\dim(\Aut(M,H))<\dim(G)$, then $\ell$ equals the
dimension of a proper graded subalgebra $\frak b=\oplus_{i=-k}^k\frak
b_i$ of $\fg$.

\noindent
(3) If $(M,H)$ is not locally isomorphic to the canonical distribution
on $G/P$ induced by $\fg^{-1}/\frak p\subset\fg/\frak p$, then the
graded Lie subalgebra in (2) has the additional property that there is
a nonzero element in $H^2(\fg_-,\frak g)$ which is annihilated by all
elements of $\frak b_0$.
\end{kor*}
\begin{proof}
  In view of the equivalence of categories proved in theorem \ref{2.4}
  parts (1) and (2) follow directly from parts (1) and (3) of theorem
  \ref{2.5}.

For part (3), we can choose our point $u\in\Cal G$ in such a way that
$\ka_h(u)\neq 0$. Then the value in $u$ of the corresponding
equivariant function determines an element of
$H^2(\fg_-,\fg)$. Elements of $\frak a^0$ correspond to infinitesimal
automorphisms $\xi$ such that $\om(\xi)(u)\in\fg^0=\frak p$, i.e.~such
that $\xi(u)$ is vertical. Differentiating an equivariant function
along $\xi$, the result in $u$ therefore coincides with the algebraic
action of $\om(\xi)(u)$ on the value of the function. By complete
reducibility of $H^2(\fg_-,\fg)$, this action depends only on the
component in $\fg_0=\fg^0/\fg^1$. 
\end{proof}

\subsection*{Remark \thesubsection} (1) By definition any infinitesimal
automorphism $\xi$ of a Cartan geometry of type $(G,P)$ is a right
invariant vector field and hence projects to an vector field
$\underline{\xi}$ on $M$. For geometries equivalent to distributions,
this is the ``actual'' infinitesimal automorphism, i.e.~its local
flows preserves the distribution. The filtration from part (3) of the
theorem has a nice interpretation in terms of $\underline{\xi}$.
Namely, for $u\in\Cal G$ and $x=p(u)\in M$, we have
$\om(\xi)(u)\in\fg^i$ if and only if $\underline{\xi}(u)\in T^i_xM$
for all $i<0$. Likewise $\om(\xi)(u)\in\fg^0$ if and only if
$\underline{\xi}(u)=0$, so we have a fix point for the one--parameter
group of automorphisms generated by $\xi$. The fact that
$\om(\xi)(u)\in\fg^i$ for some $i>0$ can be interpreted (in a certain
sense) as higher order vanishing of $\underline{\xi}$ in $x$.

\noindent
(2) Part (3) of Corollary \ref{2.5} shows that the dimension of $\frak
b_0$ is bounded by the largest possible dimension of an annihilator of
a nonzero element in $H^2(\fg_-,\fg)$, which can be easily determined
using representation theory as follows. The Lie algebra $\fg_0$ is
always reductive, and in the cases studied in this paper its center is
one--dimensional and acts non--trivially on $H^2(\fg_-,\fg)$.

For a complex semisimple Lie group $H$ and a complex irreducible
representation $V$, one can consider the action of $H$ on the
projectivization $\Bbb P(V)$. It is well known (see \cite[chapter
23]{Fulton-Harris}) that the orbit of the line through a highest
weight vector is the unique $H$--orbit of smallest dimension, and the
stabilizer of this line is a parabolic subgroup of $H$, whose type can
be read off the highest weight of $V$. Consequently, the stabilizer in
$H$ of a highest weight vector in $V$ has the largest possible
dimension among all stabilizers of non--zero vectors. Since the
parabolic subgroup acts non--trivially on the highest weight line,
this stabilizer has codimension one in the parabolic subgroup.

If $H$ is reductive with one--dimensional center which acts
non--trivially (and by Schur's lemma by a scalar) on $V$, then one can
look at the corresponding parabolic subgroup in the semisimple part of
$H$. Since both this parabolic subgroup and the center act
non--trivially on the highest weight line, it follows that the
stabilizer of a highest weight vector in $H$ has the same dimension as
the parabolic subgroup in the semisimple part. 

To apply this to our situation, we only have to notice that via
complexification, the complex dimension of the stabilizer of a
non--zero element gives an upper bound for the real dimension of the
stabilizer of a non--zero element. Hence to obtain the upper bounds
for $\dim(\frak b_0)$ one only has to compute the dimensions of
parabolic subalgebras in the complexification of the semisimple part
of $\fg_0$, whose type can be read off the highest weight of
$H^2(\fg_-,\fg)$, which is the output of Kostant's theorem. 

\section{Examples}\label{3} 

\subsection{Generic rank $n$ distributions in dimension
  $\frac{n(n+1)}{2}$}\label{3.1}
For $n\geq 3$ consider $\Bbb R^{2n+1}$ with the (indefinite) inner product 
$$
\langle v,w\rangle=v_{n+1}w_{n+1}+\sum_{i=1}^{n}v_iw_{n+1+i}+
\sum_{i=1}^{n}v_{n+1+i}w_{i}. 
$$
By definition, the subspaces generated by the first $n$ respectively
by the last $n$ vectors in the standard basis are isotropic, which
shows that the inner product has split signature $(n+1,n)$. The
orthogonal Lie algebra $\fg\cong\frak{so}(n+1,n)$ for this inner
product has the form  
$$
\fg=\left\{\begin{pmatrix}A & v & B\\ w & 0 & -v^t\\
C&-w^t&-A^t\\\end{pmatrix}:
\begin{array}{l}
A\in\frak{gl}(n,\Bbb R), C,B\in\mathfrak{o}(n)\\
v\in\Bbb R^n, w\in\Bbb R^{n*}
\end{array}\right\}.
$$
This admits an obvious grading of the form 
$$
\begin{pmatrix}\g_0 & \g_1 & \g_2\\\g_{-1} &\g_0 &\g_1 \\
\g_{-2}&\g_{-1} &\g_0\\\end{pmatrix}, 
$$
with blocks of sizes $n$, $1$, and $n$. The associated parabolic
subalgebra $\frak p=\fg^0$ is the stabilizer of the $n$--dimensional
isotropic subspace generated by the first $n$ basis vectors. From this
representation it is evident that $\fg_0\cong\frak{gl}(n,\Bbb R)$. The
adjoint action makes each $\fg_i$ into a $\fg_0$--module and the
bracket on $\fg$ induces homomorphisms of $\fg_0$--modules. From the
matrix representation it is evident, that  $\fg_{1}\cong\Bbb R^n$,
$\fg_{-1}\cong\Bbb R^{n*}$ as $\fg_0$--modules. Further, the bracket
induces isomorphisms $\La^2\fg_{\pm 1}\to\fg_{\pm 2}$ and
$\fg_{-1}\otimes\fg_1\to\fg_0$. Finally, the restrictions $\fg_{\pm
  1}\x\fg_{\mp 2}\to\fg_{\mp 1}$ of the bracket are induced by
$(w,B)\mapsto wB$ and $(v,C)\mapsto Cv$, respectively. 

The general tools mentioned in \ref{2.1} show that $H^1(\fg_-,\fg)$ is
concentrated in negative homogeneities in this case, see also
\cite{katja}.  Putting $G=\Aut(\fg)$ and $P=\Aut_f(\fg)$ we can thus
apply the results from \ref{2.4}. The corresponding bracket generating
distributions are rank $n$ distributions on manifolds of dimension
$n(n+1)/2$. If $H\subset TM$ is such a distribution, then the
condition that the bundle of symbol algebras is locally trivial and
modelled on $\fg_-$ simply means that for each $x\in M$, the Levi
bracket induces an isomorphism $\La^2H_x\to T_xM/H_x$. Since we know
that one such distribution exists (on $G/P$, which can be viewed as
the Grassmannian of all isotropic $n$--dimensional subspaces in $\Bbb
R^{2n+1}$), this is evidently a generic condition. This clearly is the
only generic type of rank $n$ distributions in dimension
$\frac{n(n+1)}{2}$. For the minimal value $n=3$, we obtain generic
rank 3 distributions in dimension 6. These have been studied by
R.~Bryant in his thesis, see \cite{Bryant1,Bryant2}.

\begin{thm*}
  Let $M$ be a smooth manifold of dimension $n(n+1)/2$, $H\subset TM$
  a generic distribution of rank $n$, and let $\Aut(M,H)$ be the group
  of all diffeomorphisms of $M$ which are compatible with $H$. 

\noindent
(1) $\Aut(M,H)$ is a Lie group of dimension at most $2n^2+n$.

\noindent
(2) If $\dim(\Aut (M,H))<2n^2+n$, then $\dim(\Aut (M,H))\leq2n^2-n+1$,
i.e. the dimension has to drop by at least $2n-1$. 

\noindent
(3) If $(M,H)$ is not locally isomorphic to the canonical
distribution on $G/P$, then 
$$\dim(\Aut(M,H))\leq
\left\{
\begin{array}{cc}
13 & n=3\\
27 & n=4\\
2n^2-3n+6 & 5\leq n\leq 8\\
2n^2-3n+5 & n\geq 9
\end{array}\right.$$ 
\end{thm*}
\begin{proof}
(1) follows immediately from part (1) of Corollary \ref{2.5}. In view
of part (2) of that Corollary, we can prove (2) by showing that any
proper graded subalgebra $\frak b=\oplus_{i=-2}^2\frak b_i$ of $\fg$
has dimension at most $2n^2-n+1$. Suppose that $\frak b$ is such a
subalgebra and put $d_j=\dim(\frak b_j)$ for all $j=-2,\dots,2$. 

Let us first assume that $d_{-1}=n$, i.e.~$\frak b_{-1}=\frak g_{-1}$.
Then $\fg_{-2}=[\fg_{-1},\fg_{-1}]\subset\frak b$, so all of $\fg_-$
is contained in $\frak b$. Then we must have $\frak b_1\neq\fg_1$,
since otherwise also $\fg_2=[\fg_1,\fg_1]$ and
$\fg_0=[\fg_{-1},\fg_1]$ would be contained in $\frak g$.

Hence $\ell:=d_1<n$. The group $G_0=\Aut_{gr}(\fg_-)$ is isomorphic to
$GL(\fg_{-1})$ and conjugating $\frak b$ by an appropriate element of
this group, we may assume that $\frak b_1=\Bbb R^\ell\subset\Bbb
R^n=\fg_1$.  Now $[\frak b_{-1},\frak b_2]\subset\frak b_1$ and since
$\frak b_{-1}=\fg_{-1}$ the description of the bracket shows that the
matrices in $\frak b_2$ may have nonzero entries only in the first
$\ell$ rows, so $d_2\leq \ell(\ell-1)/2$ by skew symmetry. In
particular, if $d_1=0$ then $d_2=0$ and we have already lost
$n(n+1)/2\geq 2n$ dimensions. So we may assume that $1\leq\ell\leq
n-1$. Then $\frak b_0\subset\fg_0=\frak{gl}(n,\Bbb R)$ must stabilize
$\frak b_1$, so $d_0\leq n^2-(n-\ell)\ell$, and this becomes
maximal for $\ell=1$ and $\ell=n-1$, whence $d_0\leq n^2-n+1$. From
above, we know that $d_2\leq\frac{(n-1)(n-2)}{2}=\frac{n^2-3n+2}{2}$,
so we conclude that 
$$
\dim(\frak b)\leq\frac{n(n-1)}{2}+n+(n^2-(n-\ell)\ell)+i+
\frac{\ell(\ell-1)}{2}\leq2n^2-n+1,
$$
with equality attained for $\ell=n-1$. 

If $d_1=n$, the situation is completely symmetric and we get the same
bound on $\dim(\frak b)$, so it remains to consider the case that both
$d_1$ and $d_{-1}$ are smaller than $n$. We claim, that in these cases
already $d_{-1}+d_0+d_1\leq
\dim(\fg_{-1})+\dim(\fg_0)+\dim(\fg_1)-2n+1$. This is evidently true
if $d_{-1}=d_1=0$, so by symmetry we may assume that for $d_1=\ell$ we
have $0<\ell<n$. Then we know that $d_0\leq n^2-(n-\ell)\ell$, so we
have lost at least $n-1$ dimensions already. If $d_{-1}\leq n-\ell$
then we have also lost at least $n$ dimensions from $d_{\pm 1}$ and
we are done again. 

Hence we are left with the case that $\ell':=d_{-1}>n-\ell$. Since
$\fg_{-1}$ and $\fg_1$ are dual $\fg_0$--modules, the annihilator of
$\frak b_{-1}$ in $\frak g_1$ is a subspace of dimension
$n-\ell'<\ell$ which by construction must be invariant under $\frak
b_0$. Hence $\frak b_0$ has to preserve two subspaces of different
dimensions. Preserving the $\ell$--dimensional subspace forces
$d_0\leq \dim(\fg_0)-(n-1)$. If the second subspace is contained in
the first one, we loose $\ell-1$ more dimensions, so the total loss
adds up to $(n-1)+\ell-1+n-\ell+n-\ell'=3n-2-\ell'$, and since
$\ell'<n$, this is at least $2n-1$. If the two subspaces are not
nested but have non--trivial intersection, then the same argument
applies to one of them and the intersection. Finally, if the two
subspaces are transverse, then each of them causes a reduction of
$(n-1)$ dimensions for $\frak b_0$.

\noindent
(3) Denote by $\frak b=\oplus_{i=-2}^2\frak b_i$ the graded subalgebra
of $\frak g$ associated to the Lie algebra of the automorphism group,
cf. Corollary \ref{2.5}.  From part (3) of this Corollary we know that
$\frak b_0$ annihilates a nonzero element in $H^2(\frak g_-,\frak g)$.

If $n=3$ the $\frak g_0$-module $H^2(\frak g_-,\frak g)$ is the
irreducible component of highest weight in $\Bbb R^3\otimes\Bbb
R^3\otimes\La^2\Bbb R^3\otimes(\Bbb R^3)^*$. Following Remark
\ref{2.5} (2), $\dim(\frak b_0)$ is bounded by the dimension of the
Borel subalgebra of the complexification of $\fg_0$, so $\dim(\frak
b_0)\leq 5$.

Suppose first that $d_{-1}=3$. Since $\g_0\cong\g_{-1}\otimes\g_1$ and
$d_0\leq5$, we obtain $d_1\leq1$.  Hence $d_2=0$ and $\dim(\frak
b)\leq12$. The same argument applies to $d_1=3$. If $d_{-1}\leq2$,
then $\g_0\cong\g_{-1}\otimes\g_1$ and $d_0\leq 5$ imply $d_1\leq2$.
If $d_{-1}=2$, we may assume that $\frak b_{-1}=\Bbb R^2\subset\Bbb
R^3$.  Since $d_1\leq2$ and $\frak g_2\times\frak
b_{-1}\rightarrow\frak g_1$ is surjective, we must have $d_2\leq2$.
For $d_1\leq1$ we are done, since we already lost $8$ dimensions.  If
$d_1=2$ we conclude from $[\frak b_{-2},\frak b_1]\subset\frak b_{-1}$
that $d_{-2}\leq2$ and therefore $\dim(\frak b)\leq13$ also in this
case. For $d_{-1}\leq 1$ and $d_1\leq1$ one already has $\dim(\frak
b)\leq13$.  If $d_{-1}\leq 1$ and $d_1>1$ we can use the arguments
above to see $\dim(\frak b)\leq 13$, thus completing the proof for
$n=3$.

For $n>3$ the $\frak g_0$-module $H^2(\frak g_-,\frak g)$ is the
irreducible component of highest weight in $\Bbb R^n\otimes\La^2\Bbb
R^n\otimes\Lambda^2(\Bbb R^n)^*$. Following Remark \ref{2.5} (2), we
can compute the dimension of parabolic subalgebra in the
complexification of $\fg_0$ to conclude that $\dim(\frak b_0)\leq
n^2-4n+10$. If $d_{-1}=n$, then we have $d_1\leq\frac{n^2-4n+10}{n}$,
since $\g_0\cong\g_{-1}\otimes\g_1$.  Further we know that
$d_2\leq\frac{(d_1-1)d_1}{2}$. Putting this together we conclude
$\dim(\frak b)\leq 2n^2-7n+26-\frac{35}{n}+\frac{100}{n^2}$. The same
argument applies if $d_1=n$, so it remains to consider the case
$d_{-1},d_1\leq n-1$.  

If $0\leq d_{-1}\leq n-1$ and $2\leq d_1\leq n-1$, then the
surjectivity of $\fg_{-2}\x\frak b_1\rightarrow \fg_{-1}$ forces
$d_{-2}<\dim(\fg_{-2})$. A straightforward analysis of the possible
cases using $d_1\cdot d_{-1}\leq n^2-4n+10$ shows that we get
$\dim(\frak b)\leq 2n^2-3n+5+\frac{7}{n-1}$.  By the symmetry of the
grading one obtains the same bound for $d_{-1}\geq2$. Finally, if
$d_{-1},d_1\leq1$, we already have $\dim(\frak b)\leq 2n^2-5n+12$.

Comparing the three bounds obtained so far, we see that $\dim(\frak
b)\leq 2n^2-3n+5+\frac{7}{n-1}$ is always valid for $n\geq 5$ while
$\dim(\frak b)\leq27$ for $n=4$.
\end{proof}

\subsection*{Remark \thesubsection}
(1) From the proof we see that for each $k=1,\dots, n-1$ 
$$
\frak b^k:=\La^2(\Bbb R^n)^*\oplus(\Bbb R^n)^*\oplus
\left(\begin{smallmatrix}*&*\\0&*\\\end{smallmatrix}\right)\oplus
\Bbb R^{k}\oplus\La^2\Bbb R^{k}
$$
is a graded subalgebra of $\fg$. The algebra $\frak b^{n-1}$ is a graded
subalgebra of the maximal possible dimension $2n^2-n+1$. It turns out
that, up to conjugation, it is the unique graded subalgebra of this
dimension. For each $k$, one can actually realize $\frak b^k$ as the
Lie algebra of the automorphism group of a generic distribution.
Namely, consider the canonical distribution on the homogeneous model
$G/P$. Then $G/P$ is the space of all maximal isotropic subspaces of
$\Bbb R^{n+1,n}$. For each $k=1,\dots, n$ let $W_k$ be the isotropic
subspace spanned by the last $(n-k)$ vectors in the standard basis of
$\Bbb R^{n+1,n}$. Then the maximal isotropic subspaces which intersect
$W_k$ only in $0$ form an open subset $U_k$ of $G/P$. Restricting the
canonical distribution to this subspace, Corollary \ref{2.4} shows
that the automorphism group of this restriction is the subgroup of $G$
consisting of all elements whose action on $G/P$ preserves $U_k$. It
is elementary to show that this coincides with the stabilizer of $W_k$
in $G$, and hence has Lie algebra isomorphic to $\frak b^k$. 

\noindent
(2) The fact that $G_0=GL(\fg_{-1})$ implies that sub--Riemannian
metrics on generic distributions of the type we consider have no
pointwise invariants. Hence for any sub--Riemannian structure on such
a distribution \cite{Morimoto:sub-Riemann} constructs a canonical
Cartan connection.  

\subsection{Generic rank two distributions in dimension
  five}\label{3.2} 

This is the classical example studied in Cartan's article
\cite{Cartan:five} from 1910. The simple Lie algebra in question is
the split real form of the exceptional Lie algebra of type $G_2$.
Although it is not difficult to describe an explicit matrix
representation for this Lie algebra, all the information we need can
be directly obtained from the root system of type $G_2$. There are two
simple roots, $\al_1$ and $\al_2$, and the other positive roots are
$\al_1+\al_2$, $2\al_1+\al_2$, $3\al_1+\al_2$, and $3\al_1+2\al_2$.
The grading we are interested in comes from the coefficient of the
short simple root $\al_1$. Hence this is a $|3|$--grading, and the
dimensions of the grading components are $\dim(\fg_{\pm
  3})=\dim(\fg_{\pm 1})=2$, $\dim(\fg_{\pm 2})=1$, and
$\dim(\fg_0)=4$. This decomposition works for the complex simple Lie
algebra of type $G_2$ as well as for its split real form. The root
decomposition also implies immediately that the Lie bracket induces
isomorphisms $\fg_0\cong\frak{gl}(\fg_{-1})$, $\La^2\fg_{\pm
  1}\cong\fg_{\pm 2}$, and $\fg_{\pm 1}\otimes\fg_{\pm 2}\to\fg_{\pm
  3}$. Note that together with the dimensions of the components, the
last two statements completely determine the structure of the
subalgebra $\fg_-$. Finally, for $i=1,2,3$ the components $\fg_i$ and
$\fg_{-i}$ are dual $\fg_0$--modules. Using this, we can identify
$\fg_{\pm 3}$ with $\fg_{\pm 2}\otimes\fg_{\pm 1}$ and $\fg_{\mp 1}$
with $\fg_{\pm 1}^*$, and under these identifications, the bracket
$\fg_{\pm 3}\otimes\fg_{\mp 1}\to\fg_{\pm 2}$ is induced by the dual
pairing $\fg_{\pm 1}\otimes \fg_{\pm 1}^*\to\Bbb R$.

Now suppose that $M$ is a smooth manifold of dimension five and
$H\subset TM$ is a distribution of rank $2$. Since $\La^2H$ then has
rank one, we see that for $x\in M$ the subspace spanned by sections of
$H$ and brackets of two such sections can have dimension at most
three, and forming the bracket with another section of $H$, one can
get at most two additional dimensions. The distribution $H$ is called
generic if and only if brackets of sections of $H$ of length at most
three span all of the tangent space, i.e.~if and only if $H$ has small
growth vector $(2,3,5)$. This is evidently equivalent to the fact that
the Levi bracket induces isomorphisms $\La^2H_x\to T_x^{-2}M/H_x$ and
$(T_x^{-2}M/H_x)\otimes H_x\to T_xM/T_x^{-2}M$, i.e.~that each symbol
algebra is isomorphic to $\fg_-$. Putting $G:=\Aut(\fg)$ and
$P:=\Aut_f(\fg)$, Theorem \ref{2.4} implies that regular normal
parabolic geometries of type $(G,P)$ are equivalent to generic rank
two distributions on five--manifolds. 

\begin{thm*}
  Let $M$ be a smooth manifold of dimension five, $H\subset TM$ a
  generic distribution of rank two, and $\Aut(M,H)$ the automorphism
  group of this distribution. 

\noindent
(1) $\Aut(M,H)$ is a Lie group of dimension at most $14$. 

\noindent
(2) If $\dim(\Aut(M,H))<14$, then $\dim(\Aut(M,H))\leq 9$. 

\noindent
(3) If $(M,H)$ is not locally isomorphic to the canonical distribution
on $G/P$, then $\dim(\Aut(M,H))\leq 8$.
\end{thm*}
\begin{proof}
  For (2) it remains to show that any proper graded subalgebra of
  $\fg$ has dimension at most nine. Suppose that $\frak
  b=\oplus_{i=-3}^3\frak b_i$ is such a subalgebra and put
  $d_j:=\dim(\frak b_j)$ for $j=-3,...,3$. Let us first assume that
  $\fg_{-1}\subset\frak b$, i.e.~that $d_{-1}=2$. Since $\fg_-$ is
  generated by $\fg_{-1}$, this implies that $\fg_-\subset\frak b$.
  Hence $\frak g_1$ cannot be contained in $\frak b$, so $d_1<2$.
  
  Now the bracket induces an isomorphism
  $\g_{2}\otimes\g_{-1}\to\g_{1}$, so we must have $d_2=0$. But then
  also $d_3=0$, since bracketing with a nonzero element of $\fg_3$ is
  a surjection $\fg_{-1}\to\fg_2$. If $d_1=0$, then we conclude
  $\dim(\frak b)\leq 9$. If $d_1=1$, then the fact the the adjoint
  action of $\frak b_0$ on $\fg_1$ must preserve the one dimensional
  subspace $\frak b_1\subset\fg_1$ implies that $d_0\leq 3$, and hence
  we again get $\dim(\frak b)\leq 9$, so the case $d_{-1}=2$ is
  complete. By symmetry this also applies if $d_1=2$, so we are left
  with the case that $d_{-1}, d_1\leq 1$.
  
  If $d_{-1}=0$, then either $d_2=0$ or $d_{-3}=0$, and likewise
  $d_1=0$ implies that $d_{-2}=0$ or $d_3=0$. In particular,
  $d_{-1}=d_1=0$ implies $\dim(\frak b)\leq 8$. If $d_{-1}=0$ and
  $d_1=1$, then $d_0\leq 3$ implies that $\dim(\frak b)\leq 9$. By
  symmetry, this also holds if $d_1=0$ and $d_{-1}=1$. Finally, if
  $d_{-1}=d_1=1$, then $d_0\leq 3$, either $d_{-3}=1$ or $d_2=0$, and
  either $d_3=1$ or $d_{-2}=0$, so again $\dim(\frak b)\leq9$.

\noindent
(3) Let $\frak b=\oplus_{i=-3}^3\frak b_i$ be the graded subalgebra of
$\g$ associated to the Lie algebra of the automorphism group. By part
(3) of Corollary \ref{2.5}, $\frak b_0$ stabilizes a nonzero element
in the irreducible $\frak g_0$- module $H^2(\frak g_-,\frak g)\cong
S^4(\frak g_1)$. Following Remark \ref{2.5} (2), we get $\dim(\frak
b_0)\leq2$ and the arguments from the proof of (2) show that
$\dim(\frak b)\leq 8$.
\end{proof}

The simplest example of a proper graded subalgebra of the maximal
possible dimension $9$ is given by $\fg_-\oplus\fg_0$. This can be
realized as the Lie algebra of the automorphism group of the
complement of a point in the homogeneous model $G/P$. Similarly,
$\fg_-\oplus\frak b_0\oplus\frak b_1$ for a line $\frak
b_1\subset\fg_1$ and its stabilizer $\frak b_0$ in $\frak g_0$ can be
easily realized. The homogeneous model $G/P$ can be viewed as the
space of null lines in the seven dimensions space of purely imaginary
split octonions and $\frak b$ can be realized as the Lie algebra of
the automorphism group of the complement of a fixed null line. It is
easy to give a complete description of all the graded subalgebras of
$\fg$ of dimension $9$, and all of them can be realized, see
\cite{Neusser}.

\subsection{Generic rank four distributions in dimension
  seven}\label{3.3} 

The last example we consider is rank four distributions with small
growth vector $(4,7)$ on manifolds of dimension seven. For such a
manifold, the Levi bracket in a point $x$ is a map $\La^2H_x\to
T_xM/H_x$. Choosing an isomorphism $H_x\to\Bbb R^4$ and
$T_xM/H_x\to\Bbb R^3$, we have to deal with the space $L(\La^2\Bbb
R^4,\Bbb R^3)$ of linear maps, which has dimension $18$.  Changing the
identifications of $H_x$ with $\Bbb R^4$ and of $T_xM/H_x$ with $\Bbb
R^3$ is expressed by the natural action of the group $GL(4,\Bbb R)\x
GL(3,\Bbb R)$, which has dimension $25$. Thus the Levi bracket in $x$
corresponds to a well defined orbit of the natural action of
$GL(4,\Bbb R)\x GL(3,\Bbb R)$ on $L(\La^2\Bbb R^4,\Bbb R^3)$.

In this picture, a distribution is stable (up to isomorphism) under
small perturbations if and only if the corresponding orbit is open in
$L(\La^2\Bbb R^4,\Bbb R^3)$. It turns out (see
\cite[7.12]{Montgomery}) that there are exactly two open orbits, and
hence two types of generic distributions. From the dimension count, we
see that this is equivalent to the fact that the stabilizer in
$GL(4,\Bbb R)\x GL(3,\Bbb R)$ of the corresponding map has dimension
seven. We can also view an element of $L(\La^2\Bbb R^4,\Bbb R^3)$ as
defining the structure of a graded Lie algebra of $\Bbb R^4\oplus\Bbb
R^3$, and then the stabilizer in $GL(4,\Bbb R)\x GL(3,\Bbb R)$ is
exactly the automorphism group of this graded Lie algebra. Hence
generic distributions are in bijective correspondence with graded Lie
algebras structures on $\Bbb R^4\oplus\Bbb R^3$ whose automorphism
group has dimension seven (which is the minimal possible dimension). 

To relate this to parabolic geometries, consider the complex simple
Lie algebra $\fg^{\Bbb C}:=\frak{sp}(6,\Bbb C)$ of endomorphisms of
$\Bbb C^6$ which preserve a non--degenerate skew symmetric bilinear
form. For this form, we use
$(z,w)\mapsto\sum_{j=1}^6(-1)^jz_jw_{6-j}$. For a matrix $A\in
\frak{gl}(2,\Bbb C)$ we denote by $\overline{A}\in\frak{gl}(2,\Bbb C)$
the classical adjoint, so $\overline{A}A=A\overline{A}=\det(A)\Bbb
I$. Using this, we can realize $\fg^{\Bbb C}$ as
$$
\left\{\begin{pmatrix}A_{11}&A_{12}&A_{13}\\A_{21}&A_{22}&
-\overline{A}_{12}\\A_{31}&-\overline{A}_{21}&-\overline{A}_{11}
\end{pmatrix}:
\begin{array}{l}
A_{13},A_{22},A_{31}\in \frak{sl}(2,\Bbb C)\\
A_{11},A_{12},A_{21}\in\frak{gl}(2,\Bbb C)
\end{array}
\right\}
$$
In this realization, we have an evident $|2|$--grading, with
$\fg^{\Bbb C}_0$ corresponding to the entries $A_{11}$ and $A_{22}$,
$\fg^{\Bbb C}_{-2}$ to $A_{31}$, $\fg^{\Bbb C}_{-1}$ to $A_{21}$,
$\fg^{\Bbb C}_1$ to $A_{12}$, and $\fg^{\Bbb C}_2$ to $A_{13}$. Hence
the complex dimensions of the $\pm 1$--components are four, the $\pm
2$--components have complex dimension three, while the $0$--component
has complex dimension seven. The Lie bracket of $\fg^{\Bbb C}_-$ thus
defines a map $\La^2\Bbb C^4\to\Bbb C^3$. Now the cohomological
condition from \ref{2.4} is satisfied for any real form of the given
grading on $\fg^{\Bbb C}$. From \ref{2.4} we know that if
$\fg=\oplus_{i=-2}^2\fg_i$ is such a real form and $G=\Aut(\fg)$, then
the subgroup $G_0\subset G$ with Lie algebra $\fg_0$ is isomorphic to
the automorphism group of the graded Lie algebra $\fg_-$. Since
$\dim(\fg_0)=7$, any such real form corresponds to a generic rank four
distribution in dimension seven, and by Theorem \ref{2.4} such
distributions then are equivalent to regular normal parabolic
geometries of type $(G,P)$. 

It turns out that there are two real forms of this grading. One is the
split real form $\frak{sp}(6,\Bbb R)$. With respect to the obvious
real analog of the bilinear form used above, we get the same matrix
representation as for $\fg^{\Bbb C}$ above but with the blocks lying
in $\frak{gl}(2,\Bbb R)$ respectively $\frak{sl}(2,\Bbb R)$. The other
real form is obtained by passing to quaternionic matrices. If we
identity $\Bbb C^6$ with $\Bbb H^3$, then our skew symmetric form
comes from a quaternionic Hermitian form of signature $(2,1)$. Maps
preserving both the quaternionic structure and the skew form also
preserve the quaternionic Hermitian form, and hence form a Lie algebra
isomorphism to $\frak{sp}(2,1)$. The matrix representation in this
case is exactly as above but with the quaternions $\Bbb H$ replacing
$\frak{gl}(2,\Bbb R)$, the purely imaginary quaternions $\im(\Bbb H)$
replacing $\frak{sl}(2,\Bbb R)$, and the quaternionic conjugation
instead of the classical adjoint. 

The deeper reason for this dichotomy is that up to isomorphism there
are exactly two real composition algebras of dimension four, namely
the quaternions with their standard positive definite quadratic form,
and the split--quaternions, which are isomorphic to the algebra of
$2\x 2$--matrices, with the quadratic form given by the
determinant. This quadratic form is preserved by the action of $G_0$
up to scale. Passing to distributions this means that any generic rank
four distribution in dimension seven admits a canonical conformal
class of inner products (i.e.~a canonical conformal sub--Riemannian
structure) which is either definite or of split signature. We refer to
these cases as elliptic type respectively hyperbolic type. 

\subsection{Results for elliptic type}\label{3.4}
Here we have $\fg=\frak{sp}(2,1)$, $\fg_{\pm 1}\cong\Bbb H$, and
$\fg_{\pm 2}\cong\im(\Bbb H)$. The bracket $\La^2\fg_{\pm
  1}\to\fg_{\pm 2}$ is given by $(p,q)\mapsto p\bar q- q\bar p$ for
$p,q\in\Bbb H$. The brackets $\fg_{\pm 2}\x\fg_{\mp 1}\to \fg_{\pm_1}$
is given by $(p,q)\mapsto p\bar q$. Finally, the adjoint action on
$\fg_{-1}\cong\Bbb H$ identifies $G_0$ with the conformal group
$CSO(4)$. 

\begin{thm*}
  Let $M$ be a smooth manifold of dimension seven, let $H\subset
  TM$ be a generic rank four distribution of elliptic type, and let
  $\Aut(M,H)$ be the automorphism group of $H$.

\noindent
(1) $\Aut (M,H)$ is a Lie group of dimension at most $21$. 

\noindent
(2) If $\dim(\Aut(M,H))<21$, then $\dim(\Aut(M,H))\leq 14$.

\noindent
(3) If $(M,H)$ is not locally isomorphic to the canonical distribution
on $G/P$, then $\dim(\Aut(M,H))\leq 12$.
\end{thm*}
\begin{proof}
  As before, we only have to prove (2) and (3), and for (2) we have to
  show that any proper graded subalgebra $\frak b\subset\fg$ has
  dimension at most $14$. We put $\frak b=\oplus_{j=-2}^2\frak b_j$
  and $d_j:=\dim(\frak b_j)$.
  
  Let us first assume that $\frak b_{-1}=\fg_{-1}$, i.e.~that
  $d_{-1}=4$. Then $\frak b_{-2}=\fg_{-2}$, and to obtain a proper
  subalgebra, we have to have $d_1<4$. Since the bracket with a
  non--zero element of $\fg_2$ defines a surjection
  $\fg_{-1}\to\fg_1$, we see that we must have $d_2=0$, which in turn
  implies that $d_1\leq 1$. If $d_1=0$, then we have $\dim(\frak
  b)\leq 14$. If $d_1=1$, then $\frak b_0$ has to stabilize a line in
  $\fg_1$, which forces $d_0\leq 4$, so we get $\dim(\frak b)\leq 12$.
  If $d_1=4$, then the result follows in the same way by symmetry.
  
  Let us next assume that $d_{-1}=3$. Conjugating by an element of
  $G_0$, we may assume that $\frak b_{-1}$ is spanned by $1,i,j\in\Bbb
  H$. This shows that $\frak b_{-2}=\frak g_{-2}$, and then the fact
  that $[\frak b_1,\fg_{-2}]\subset\frak b_{-1}$ forces $d_1\leq 1$.
  This is only possible if $d_2=0$, and $\dim(\frak b)\leq 14$ readily
  follows. The same argument applies if $d_1=3$.
  
  If $d_{-1}=2$, then we may assume that $\frak b_{-1}$ is spanned by
  $1$ and $i$, and then conclude that $d_{-2}\geq 1$. If $d_{-2}=1$,
  then we must have $d_1\leq 2$, and this is only possible if $d_2\leq
  1$, so we have lost eight dimensions already. If $d_{-1}>1$, then
  $d_1\leq 1$ and $d_2=0$, and again we have lost eight dimensions.
  The case $d_1=2$ can be treated in the same way.
  
  If $d_{-1}=d_1=1$, then both $d_2$ and $d_{-2}$ must be $\leq 1$, so
  we have lost 10 dimensions. If $d_{-1}=1$ and $d_1=0$, then we also
  must have $d_2=0$, and again we are done. Of course, $d_{-1}=d_1=0$
  already implies a loss of eight dimensions.

\noindent
(3) In the proof of (2) we have used restrictions on $d_0$ only in one
point, and there we obtained $\dim(\frak b)\leq 12$. Hence it suffices
to prove that in the setting of (3) we have to have $d_0\leq 5$, since
this causes a loss of two more dimensions compared to (2).  But this
follows immediately from Remark \ref{2.5} (2) since the two
irreducible components in $H^2(\frak g_-,\frak g)$ are non--trivial,
and the maximal parabolic subalgebras in the complexification of
$\fg_0$ have dimension $5$.
\end{proof}

The simplest example of a proper graded subalgebra of $\fg$ of the
maximal possible dimension $14$ is of course the parabolic subalgebra
$\frak p$. Since the group $P$ is the stabilizer of the base point $o$ in
the homogeneous model $G/P$, it is the subgroup of all elements that
preserve the open subset $G/P\setminus\{o\}$. From Corollary \ref{2.4}
we conclude that $P$ is the automorphism group of the canonical
distribution on this open subset. 

\subsection{Results for hyperbolic type}\label{3.5}
The description of the brackets is as in the elliptic case but
replacing the quaternions $\Bbb H$ by the algebra $M_2(\Bbb R)$ of
real $2\x 2$--matrices, the imaginary part by the subspace
$\frak{sl}(2,\Bbb R)$ of tracefree matrices, and the conjugation of
quaternions by $\overline{A}=\Cal CA$, the classical adjoint of $A$.
The adjoint action of $\fg_{\pm 1}\cong\Bbb R^4$ identifies $G_0$ with
the conformal group $CSO(2,2)$ of split signature $(2,2)$, and the
adjoint action on $\fg_{\pm 2}\cong\Bbb R^3$ maps $G_0$ onto
$CSO(1,2)$. 

\begin{thm*}
  Let $M$ be a smooth manifold of dimension seven, let $H\subset
  TM$ be a generic rank four distribution of hyperbolic type, and let
  $\Aut(M,H)$ be the automorphism group of $H$.

\noindent
(1) $\Aut (M,H)$ is a Lie group of dimension at most $21$. 

\noindent
(2) If $\dim(\Aut(M,H))<21$, then $\dim(\Aut(M,H))\leq 16$.

\noindent
(3) If $(M,H)$ is not locally isomorphic to the canonical distribution on $G/P$, then 
$\dim(\Aut(M,H))\leq 15$.
\end{thm*}
\begin{proof}
  To prove (2) we have to  show that any
  proper graded subalgebra $\frak b\subset\fg$ has dimension at most
  $15$. We put $\frak b=\oplus_{j=-2}^2\frak b_j$ and $d_j:=\dim(\frak
  b_j)$. 
  
  If $d_{-1}=4$, then $\fg_-\subset\frak b$ and we must have $d_1<4$.
  We have seen that the bracket $\fg_2\x\fg_{-1}\to\fg_1$ is given by
  $(A,B)\mapsto A\Cal C(B)$. Hence in contrast to the elliptic case,
  the map $\fg_{-1}\to\fg_1$ defined by the bracket with a fixed
  element of $\fg_2$ is only surjective for invertible elements. Since
  any subspace of $\frak{sl}(2,\Bbb R)$ of dimension at least two
  contains an invertible matrix, we obtain $d_2\leq 1$.  Moreover, we
  know that the bracket $\fg_1\x\fg_1\to\fg_2$ is given by
  $(A,B)\mapsto A\Cal C(B)-B\Cal C(A)$, the imaginary part of $A\Cal
  C(B)$. Therefore we obtain $d_1\leq2$, since otherwise we would have
  $\dim([\frak b_1,\frak b_1]))>1$. Now since $\frak b_0$ has to
  stabilize the subspace $\frak b_1\subset\frak g_1$, we conclude that
  $d_0\leq 6$ and hence $\dim(\frak b)\leq 16$. The case $d_1=4$ can
  be treated in the same way, so it remains to consider the case that
  $d_{\pm 1}\leq3$.
  
  If $d_{-1}=3$, then $\frak b_{-1}\subset\frak g_{-1}$ is three
  dimensional, and we have to distinguish cases according to the
  signature (up to sign) of the restriction of the inner product to
  this subspace. The possible signatures are $(2,1)$, $(2,0)$, and
  $(1,1)$. Since $G_0\cong CSO(2,2)$, any two subspaces with the same
  signature are conjugate. Using this, one verifies that in the first
  two cases, $[\frak b_{-1},\frak b_{-1}]=\frak g_{-2}$, so
  $d_{-2}=3$. Then the fact that $[\frak g_{-2},\frak b_1]\subset\frak
  b_{-1}$ implies that $d_1\leq 2$, and then $[\frak b_2,\frak
  b_{-1}]\subset\frak b_1$ shows that $d_2\leq 1$. This already shows
  that $\dim(\frak b)\leq 16$. 

  For the remaining signature $(1,1)$, we may assume that $\frak
  b_{-1}$ consists of all matrices of the form
  $\left\{\left(\begin{smallmatrix}a&b\\-a&c\end{smallmatrix}\right)\right\}$.
  Computing the bracket of two matrices of this form, we see that
  $\frak b_{-2}$ has to contain all matrices of the form
  $\left\{\left(\begin{smallmatrix}\al
        &\be\\0&-\al\end{smallmatrix}\right)\right\}$ and $d_{-2}\geq
  2$. If $d_{-2}=3$, then the results follows as above. Otherwise, the
  arguments on brackets as above show that $d_1\leq 3$ and $d_2\leq 2$
  and since $\frak b_0$ has to preserve non--trivial subspaces, we get
  $d_0\leq 6$ and hence $\dim(\frak b)\leq 16$. This complete the case
  $d_{-1}=3$ and hence also the case $d_1=3$, so we are left with the
  case $d_{-1},d_1\leq2$.  Suppose $0<d_{-1}\leq2$. Since $\frak b_0$
  has to stabilize $\frak b_{-1}$, we obtain at least $d_0\leq6$.
  Hence in any case $\dim(\frak b)\leq16$.  If $d_{-1}=0$ we are
  already done, since $d_1\leq2$ by assumption.

\noindent
(3) As in the proof of part (3) of Theorem \ref{3.4}, we obtain
$d_0\leq 5$, and going through the proof of part (2), we see that in
each case we loose at least one more dimension.
\end{proof}

\subsection*{Remark \thesubsection} 
The algebras $\fg$ from \ref{3.4} and here as well as their gradings
have the same complexification, so we are just dealing with two
different real forms of one complex object. Remarkably, we obtain
different bounds on the sizes of possible automorphism groups in the
elliptic and hyperbolic case. The maximal possible dimension $16$ of
an automorphism group in the hyperbolic can actually be realized, so
the two case are truly different. A realization is obtained as
follows. In the notation from \ref{3.3}, let $\frak
b\subset\frak{sp}(6,\Bbb R)$ be the subspace formed by all matrices
for which
$A_{11}=\left\{\left(\begin{smallmatrix}*&*\\0&*\end{smallmatrix}\right)\right\}$,
$A_{12}=\left\{\left(\begin{smallmatrix}*&*\\0&0\end{smallmatrix}\right)\right\}$,
and
$A_{13}=\left\{\left(\begin{smallmatrix}0&*\\0&0\end{smallmatrix}\right)\right\}$.
One immediately verifies that this is a subalgebra of dimension $16$.
Let $B$ be the corresponding subgroup of $Sp(6,\Bbb R)$, and consider
the $B$--orbit of $o=eP$ in $Sp(6,\Bbb R)/P$. Since $\frak b$ contains
$\fg_{-}$, this orbit is open. The homogeneous space $Sp(6,\Bbb R)/P$
admits a second interpretation, namely the space of all isotropic two
planes in the symplectic vector space $\Bbb R^6$. In this picture, $B$
is the stabilizer of a hyperplane in $\Bbb R^6$, which shows that $B$
cannot act transitively on $Sp(6,\Bbb R)/P$. Hence the $B$--orbit of
$eP$ in this space is a proper open subset, so its automorphism group
must be of dimension strictly less then the dimension of $Sp(6,\Bbb
R)$. Since on the other hand, $B$ is contained in this automorphism
group, Theorem \ref{3.5} implies that we must have an automorphism
group of dimension 16.

\end{document}